\documentclass{article}

\usepackage{amsmath,amsthm,amsfonts}

\usepackage{amssymb}

\usepackage{epsfig}

\usepackage{rotating}

\usepackage{subfigure}

\def\sl R{\mathscr{R}}

\def\vp{{\varphi}}
\def\hatx{{\hat{x}}}
\def\hatR{{\hat{R}}}
\def\hatvp{{\hat{\vp}}}
\def\P{{P}}
\def\Q{{Q}}

\def\pH{\partial H}
\def\eps{\epsilon}

\begin{document}
\centerline{\Large\bf PROBLEM REDUCTION, RENORMALIZATION, AND MEMORY}
\vskip14pt

\centerline{\bf Alexandre J.\ Chorin and Panagiotis Stinis}
\vskip12pt
\centerline{Department of Mathematics, University of California}
\centerline{and}
\centerline{Lawrence Berkeley National Laboratory}
\centerline{Berkeley, CA 94720}

\begin{abstract}
Methods for the reduction of the complexity of computational problems are
presented, as well as their connections to renormalization, scaling, and irreversible statistical mechanics. 
Several statistically stationary cases are analyzed; for time dependent 
problem averaging usually fails, and averaged equations 
must be augmented by appropriate memory and random forcing terms. 
Approximations are described and
examples are given. 
\end{abstract}

\section{Introduction}\label{intro}
There are many problems in science which are too complex for numerical
solution as they stand. Examples include turbulence and other problems
where multiple scales must be taken into account. Such problems must be
reduced to more amenable forms before one computes. In the present paper
we would like to summarize some reduction methods that 
have been developed in recent years, together with an account of what was 
learned in the process. It is obvious that the problem has not been fully solved, 
but we think that the examples and the
conclusions  reached so far are useful.

In general terms,  a reduction to a more  amenable form is a renormalization
group transformation, as in physics --- a
transformation of a problem into a more tractable form while keeping 
quantities of interest invariant. A renormalization group transformation
involves an incomplete similarity transformation (see below for definitions),
and thus a reduction method is a search for hidden
similarities. This a general feature of reduction methods,
and it will be illustrated in the examples. A successful problem reduction
produces a new problem which must in some asymptotic sense be similar to the
original problem. For general backgound on renormalization, see e.g.\cite{benfatto,fisher,stanley}.

In problems with strong time dependence, reduction methods resemble methods
for the analysis of thermodynamic systems not in equilibrium; indeed, those aspects
of the problem that are ignored in a reduced description conspire to 
destroy order and increase entropy.  Problem reduction for time-dependent
problems is basically renormalization group theory for non-equilibrium 
statistical mechanics. For background on such theory, see e.g. \cite{barenblatt,goldenfeld1,kevrekidis2}.

The content of the paper is as follows: In section \ref{ave} we consider Hamiltonian systems 
and their conditional expectations. 
In section \ref{block} we narrow the discussion to statistically stationary Hamiltonian systems
and recover Kadanoff real-space renormalization
groups and an interesting block Monte-Carlo method. In section \ref{kdv} we display an example that
exhibits and also extends the main features of this analysis in simple form.

In section \ref{mz}  we explain the Mori-Zwanzig formalism for the reduction of statistically
time-dependent problems. The analysis shows that averaging the equations is in general
not enough; one must take into account noise and a temporal memory. The Mori-Zwazig
formalism is rather dense, and in the sections that follow we present various special
cases in which it can be simplified, in particular when the memory is very short or very long.

For the sake of readability, we remind the reader of the rudiments of similarity
theory
\cite{barenblatt}.
Suppose a variable $a$ is a function of variables
$a_1,a_2,\ldots, a_m$,
$b_1,b_2,\ldots, b_k$, where
$a_1,\ldots,  a_m$ have independent units, for example units of length and mass,
while the units  of
$b_1,\ldots, b_k$, can be formed from the units of 
$a_1,a_2,\ldots, a_m$.
Then there exist dimensionless variables
$\Pi=\frac{a}{a_1^{\alpha_1}\cdots a_m^{\alpha_m}}$,
$\Pi_i=\frac{b_i}{a_1^{\alpha_{i1}}\cdots a_m^{\alpha_{im}}}$,
$i=1,\ldots,k$, where the $\alpha_i,\alpha_{ij}$ are simple fractions, such that
$\Pi$ is a function of the $\Pi_{i}$:
\begin{equation}
\Pi=\Phi(\Pi_1,\ldots,\Pi_k).
\end{equation}
This is just a consequence  of the requirement that a physical relationship
be independent of the size of the units of measurement. At this stage nothing can be said about the
function $\Phi$. 
Now
suppose the variables $\Pi_i$ are small or large, and assume that the
function $\Phi$ 
has a non-zero finite limit as
its arguments tend to zero or to infinity; then $\Pi\sim$ constant, and one finds a power
monomial relation between $a$ and the $a_i$.
This is a complete similarity relation. 
If the function $\Phi$ does not have the assumed limit, it may
happen that for $\Pi_1$ small or large, $\Phi(\Pi_1)=\Pi_1^{\alpha}\Phi_1(\Pi_1)+\ldots$,
where the dots denote lower order terms, $\alpha$ is a constant, the other arguments of 
$\Phi$ have been omitted and $\Phi_1$ has a finite non-zero limit. 
One can then obtain a scaling (power monomial) expression for $a$ in terms of  the $a_i$ and $b_i$,
with undetermined powers which must be found by means other than dimensional analysis. 
The resulting power relation is an {\it incomplete}  similarity relation. 
Of course one may well have functions $\Phi$ with neither kind of similarity.

Incomplete similarity expresses what is invariant under a renormalization
group; all renormalization group transformations involve incomplete similarity, 
see  the books already cited as well as \cite{benettin} written before the
notion of incomplete similarity was formalized. The exponent $\alpha$ is called
an anomalous exponent.

The paper \cite{givon1} is a survey of reduction methods organized along different
lines and can be profitably read in tandem with the present paper.

\section{Averaging a Hamiltonian system} \label{ave} 
We begin by examining what happens when one tries to reduce the complexity of
a Hamiltonian system by averaging (see also \cite{CKK1,CKL,seibold}). 
Consider a system of nonlinear ordinary differential equations,
\begin{eqnarray}
\frac{d}{dt}\varphi(t)& =&R (\varphi(t)),\nonumber\\ 
\varphi(0)&=&x,
\label{eq:system}
\end{eqnarray}
where $\varphi$ and $x$ are $n$-dimensional vectors with components
$\varphi_i$ and $x_i$, and $R$ is a vector-valued function with components
$R_i$; $t$ is time. 
To each initial value $x$ in (\ref{eq:system}) corresponds a
trajectory $\varphi(t)=\varphi(x,t)$.

Suppose that we only want to find 
$m$ of the $n$ components
of the solution vector $\varphi(t)$ without finding the $n-m$ others. 
One has to assume something about the variables that are not evaluated, and we assume
that at time t=0 we have a a joint probability density $F(x)$ for all the variables.
The variables we keep will have definite initial values $x_1,x_2,\dots,x_m$, and the 
rest of variables will then have a conditional probability density $f_m=f(x_1,\dots,x_m,x_{m+1},\dots)/Z_m$,
where $Z_m=\int_{-\infty}^{+\infty}f(x_1,\dots,x_m,x_{m+1},\dots)dx_{m+1}dx_{m+2}\cdots$ is a normalization
constant. Without some assumption about the missing variables the problem is meaningless;
this particular assumption is reasonable because in practice $f$ can often be estimated
from previous experience or from general considerations of statistical mechanics. 
The question is how to use this prior
knowledge in the evaluation of $\varphi(t)$.

Partition the vector $x$ so that $\hatx=(x_1,x_2,\dots,x_m)$, $\tilde x=(x_{m+1},\dots,
x_n)$ and $x=(\hat x,\tilde x)$, and similarly $\varphi=(\hat\varphi,\tilde\varphi),
R=(\hat  R,\tilde R)$.  In general the first $m$ components of $R$ depend on all the components of $\vp$,  $\hatR=\hatR(\varphi)=\hatR(\hat\varphi,\tilde\varphi)$;
if they do not we have a system of $m$ equations in $m$ variables and 
nothing further needs to be done. 
We want to calculate only 
the variables $\hat\varphi$; then
$(d/dt)\hat\varphi(t) =\hat R (\varphi(t))$ where the right hand
side depends on the variables $\tilde\varphi$ which are unknown at time $t$. 
We shall call the variables $\hat\varphi$ the ``resolved variables" and the 
remaining variables $\tilde\varphi$ the ``unresolved variables".

Consider in particular a Hamiltonian system as in \cite{CKK1},\cite{CKL}. There exists then
a Hamiltonian function $H=H(\varphi)$ such that for $i$ odd $R_i$, the $i$-th component
of the vector $R$ in (\ref{eq:system}) satisfies $R_i=\partial H\bigl/{\partial \varphi_{i+1}}$
while for $i$ even one has $R_i=-{\pH}\bigl/{\partial \varphi_{i-1}}$, with $n$, the size of the system, even. Assume furthermore that $f$, the initial probability
density, is  $f(\varphi)=Z^{-1}\exp(-H/T)$ 
where $T$ is a parameter, known in physics as the ``temperature", which will be set equal to one in much, but not all, of the discussion below.
In physics this density appears naturally and is known as the ``canonical" density; 
the normalizing constant $Z=Z(T)$ is the ``partition function".
This density $f$ is invariant, i.e.
sampling it and evolving the system in time commute.

A numerical analyst who wants to approximate the solution of an equation usually
starts by approximating the equation. 
If one solves for the resolved variables one has values for the variables $\hat\varphi$ available
at each instant $t$ and the best approximation should be a function of these variables; it is natural to seek a best approximation in the
mean square sense with respect to the invariant density $f$ at each time; the best approximation
in this sense is the conditional expectation 
$E[R(\varphi)|\hat\varphi]=\int e^{-H}d\tilde\varphi\bigl/\int e^{-H}d\tilde\vp$ (note that we set $T=1$). 
This conditional expectation is the orthogonal projection of $R$ onto the space of
functions of $\hat{\varphi}$ with respect to the inner product $(u,v)=E[uv]=\int u(\varphi)v(\varphi)f(\varphi)d\varphi$, where $d\varphi$ denotes integration over all the components
of $\varphi$. 
We then try to approximate the system ($\ref{eq:system}$) by:
\begin{eqnarray}
\frac{d}{dt}\hat\varphi(t) &=& E[R (\varphi(t))|\hat\varphi(t)],\nonumber \\ 
\hat\varphi(0)&=&\hat x.
\label{foop}
\end{eqnarray}

We have shown in \cite{CHK3,CKK1,CHK} that: 
(i) The new system (\ref{foop}) is also Hamiltonian:
\begin{equation}
E \left[\frac{\pH}{\partial \varphi_i}|\hat\varphi(t)\right]=\int\frac{\pH}{\partial \varphi_i}\exp(-H)d\tilde\varphi
\bigl/ 
\int \exp(-H)d\tilde\varphi
=\frac{\partial \hat H}{\partial\varphi_i},
\label{hald1}
\end{equation}
where $i\le m=$ the dimension of $\hat\varphi$, and 
\begin{equation}
\hat H=-log\int \exp(-H)d\tilde \vp
\label{renham}
\end{equation}
is the new Hamiltonian.

(ii) The new canonical density $\hat f=Z^{-1}\exp(-\hat H)$ is invariant in the 
evolution of the
new, reduced, system.

(iii)
When the data are sampled from the canonical distribution,  the distribution of $\hat \varphi$ in the new system is its
marginal distribution in the old system; equivalently, 
the partition function $Z$  is the same for the old system and  for the new system.

Now the question is, what does the solution $\hat\varphi(t)$ of (\ref{foop}) represent ? 
It does not approximate the first $m$ components of the solution
$\varphi(t)$ of (\ref{eq:system})- the components of $\hat{\varphi}$ and 
the components of $\varphi$ live in spaces of different dimension and in general
the components of the latter in those higher $n-m$ dimensions are not small.
One could hope that what the solution of (\ref{foop}) approximates is the vector 
$E[\hat\varphi(t)|\hatx]$, the best estimate of the first components of the
solution at time $t$ given the partial initial information $\hatx$. 
This is the case for linear systems (where averaging and time integration
commute), and is approximately the case for limited time in some
other special situations- nearly linear systems, some systems where the 
``unresolved variables" are fast.
However, in general this is not the case. We shall see below that a
reduced description of the solution of nonlinear systems in time
requires in general ``noise" and a ``memory".

The lack of convergence can be  understood by the following physics
argument. 
In physics a system in which the values of all
the variables are drawn from a canonical distribution is a system in thermal 
equilibrium. 
The assignment of 
definite values $\hat x$ to the variables $\hat \varphi$ at time $t=0$
amounts to taking the system out of equilibrium at $t=0$;
if the system is ergodic it will then decay to equilibrium in time, so that
all the variables become randomized and acquire the joint density $f$.
Thus 
the  predictive value of the partial initial data $\hat x$
decreases in time; all averages of the $\hat \varphi$ approach
equilibrium averages. However, the reduced system (\ref{foop}) is Hamiltonian, and the 
solutions it produces 
oscillate forever.

In Figure 1 we consider the Hald Hamiltonian system (\cite{CHK3}) with 
\begin{equation}
H=\frac{1}{2}\left(\varphi_1^2+\varphi_2^2 +\varphi_3^2+\varphi_4^2
+\varphi_1^2\varphi_3^2\right) 
\label{haldmodel}
\end{equation}
(physically, two linear oscillators with a nonlinear coupling).
We assume that $\varphi_1(0), \varphi_2(0)$ are given and sample the two other
initial data from the canonical distribution with $T=1$.

\begin{figure}
\centering
\epsfig{file=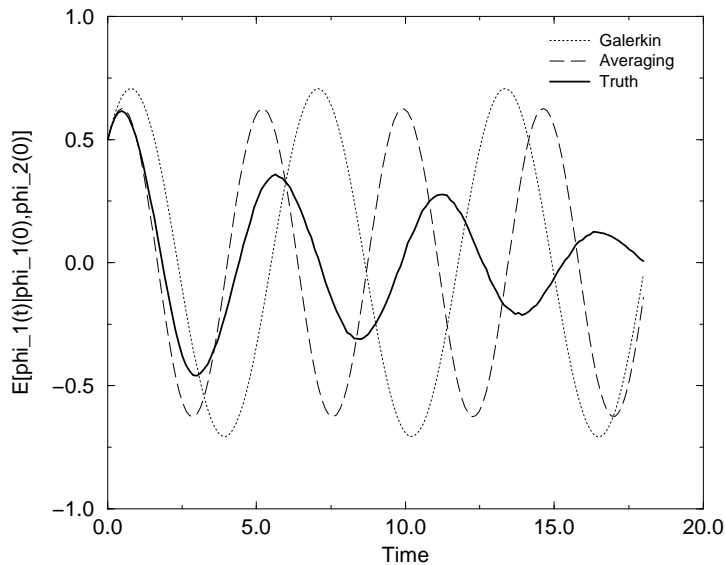,height=3in}
\caption{Comparison of the evolution of $E[\phi_1(t)|\phi_1(0),\phi_2(0)]$ (truth), to the 
prediction by the "Galerkin" approximation and the prediction by the averaging procedure described 
in the text.}
\label{fig:unkno1}
\end{figure}

In Figure 1 are displayed (1) The result  for $\varphi_1$ of a ``Galerkin" calculation in which
the unresolved variables are set to zero (this is what is implicitly done in many
unresolved computations); (2) the result of the averaging procedure just described,
and (3) the true $E[\varphi_1(t)|\hatx]$, calculated by repeatedly sampling the initial data,
solving the full system, and averaging. As one can see, averaging is initially better than
the null ``Galerkin" method, but in the long run the truth decays but the solution of the
averaged system oscillates for ever. For more detail, see \cite{CHK3}.

The procedure we have just described resembles sufficiently the averaging
methods used in some areas of engineering, for example the large-eddy
simulation methods in turbulence (see e.g. \cite{moser}) and in some multiscale problems (see e.g. \cite{kevrekidis}), to cast a very serious doubt on the broad validity of the latter.
For a description of special cases, with small fluctuations and particular structures, where this
procedure is legitimate, see \cite{givon1}.

\section{Prediction with no data and block Monte-Carlo}\label{block}

There is however a case where the construction of the preceding section
can be very useful-- when $m=0$, i.e., when one tries to predict the future 
with no initial information. Equations (\ref{eq:system}) then sample the canonical
distribution and the reduced system samples a subset of variables
without sampling the others, and, as we have seen, keep the statistics of the resolved
variables unchanged (see \cite{seibold} for an application to molecular dynamics).

To see what is happening, suppose the variables $\varphi_i$ are associated
with nodes on a regular lattice, for example, they may represent
spins in a solid, or originate in the spatial discretization
of a partial differential equation.

Divide the lattice into blocks of some fixed shape (for example, divide
a regular one-dimensional lattice into groups of two contiguous nodes). 
We had not yet specified how the variables are to be divided into
resolved and unresolved. Now decide to ``resolve" one variable per block,
and leave the others in the same block unresolved. The transformation between the old variables and the smaller set of resolved variables is a Kadanoff
renormalization group transformation \cite{kadanoff}; the Hamiltonian $\hat H$ defined
above in equation (\ref{renham}) is the renormalized Hamiltonian. We will now explain what this means.

Suppose the system described by the Hamiltonian is translation invariant.  The equations of
motion for any at any one point, say at the location labeled by $1$, have the same form as the equations of motion at any
any other point. The relation between the right hand side of the reduced system and the 
right hand side of the old system can be rewritten as:
\begin{equation}
\frac{\partial \hat H}{\partial \varphi_1}=E[\frac{\partial H}{\partial \varphi_1}|\hat\varphi],
\label{start}
\end{equation}
where the expected value is with respect to the invariant density as before. This relation is the starting
point for the actual evaluation of $\hat H$.

Hamiltonians are functions of the variables $\varphi$. They can be expanded in the form:
\begin{equation}
H=\sum_ja_j\psi_j,
\label{expandH}
\end{equation}
where the $\psi_j$ are ``elementary Hamiltonians". In a translation invariant system, where
each equation has the same form as any other, the Hamiltonian is made up of sums over $i$ of terms
of the form $h(\varphi_j\varphi_j)$ for various values of $j$, where $h$ is some function; these terms 
represent ``couplings" between variables $j$ apart;  one can then choose the elementary Hamiltonians to be
polynomials in $x_ix_{i+j}$ with a fixed $j$ in each $\psi_j$, i.e., one segregates the couplings between
variables $j$ apart into separate terms.

In a homogeneous system where there is only one variable per site it is enough to satisfy (\ref{start}) for one variable, say for $\varphi_1$. Define
$\psi'=\frac{\partial}{\partial\varphi_1}\psi$, noting that though $\psi$ is necessarily a function with
at least
as many arguments as there are components on $\varphi$, $\psi'$ can be sparse. Equation (\ref{start}) reduces to
\begin{equation}
\frac{\partial \hat H}{\partial \varphi_1}=\sum_ja_jP\psi'_j(\varphi)=\sum_j \hat{a}_j \psi'_j
(\hat{\varphi}),
\label{more1}
\end{equation}
with the projection $P$ defined as before by $Pg(\varphi)=E[g|\hat\varphi]$ for any function $g$ of $\varphi$.
Now we're almost done. One can pick a basis in $\hat L_2$, the subspace of square integrable functions that depend only
on the variables $\hat\varphi$, which consists of a subset of the set of functions $\psi'$. The right-hand
of equation (\ref{more1}) is then a linear combination of $\psi's$; integration with respect to $\varphi_1$
requires only the erasure of the primes and yields a series for $\hat H$. The elements of $\tilde\varphi$ are
now gone, and one can relabel the remaining variables $\hat\varphi$ so that the terms in the series
have exactly the same form as before; the calculation can then be repeated, yielding a sequence of 
Hamiltonians with ever fewer variables: $H, H^{(1)}=\hat H$, $H^{(2)}=\hat H^{(1)}, \dots$. The corresponding
densities $f^{n}=Z^{-1}\exp(-H^{(n)}/T)$ can in principle be sampled by any sampling scheme, for example by Metropolis sampling 
(but there are caveats, see e.g. \cite{chorin9}).

At this point we have reduced the number of variables by a factor $L$ equal to the number of
variables in each  block, but this may well seem to be a pyrrhic victory. The Hamiltonians
one usually encounters are simple, in the sense that they involve few couplings- finite
differences typically link a few neighboring variables, and so do the usual spin Hamiltonians
in physics. As one reduces the number of variables, the new Hamiltonians become more complex,
with more terms in the series (\ref{expandH}); the cost per time step of solving the equations in time or
of  the cost per move in a Metropolis sampling typically increases fast as well. To see what has
been gained one must turn to the physics literature (see e.g. \cite{kadanoff}.\cite{hohenberg}).

Consider the spatial correlation length $\ell$  which measures the range of values of 
$|j|$ over which the spatial covariances $E[\varphi_i\varphi_{i+j}]$ are non negligible,
and the correlation time $\tau$ for which the temporal covariances $E[\varphi_i(t)\varphi(t+s)]$
are non-negligible. For very large and very small values of the temperature $T$ (the variance 
parameter in the density  $f$) both the correlation time and the correlation length are small;
the properties of the system can then be found from calculations with a small number of variables and
it is not urgent to reduce the number of variables. There is a range of intermediate values of
$T$ for which the correlation length and time for are large and then the reduction is worthwhile.
There often is a value $T_c$ of $T$, the ``critical value", for which $\ell=\infty$. Values of $T$
around $T_c$ are often of great interest.

Now we can see what the reduction can accomplish. If one tries to compute averages with $T$ near
$T_c$ one finds that the cost of computation is proportional to $\tau$- one has to compute long
enough to obtain independent samples of $\varphi$, and a new independent sample will not
appear until a time $\sim\tau$ has passed. The reductions above produce a system
with smaller $\ell$ and $\tau$ and therefore computation takes less time.
Though we started with the declared goal of reducing the number of variables, what has been
produced is more interesting: a new system with shorter correlations which is more amenable to
computation. It is not the raw number of variables that matters.

The renormalization can be used with a multigrid scheme, in which one runs  up and down on different levels
of renormalization, on the finer ones to achieve accuracy and  the cruder ones to move fast from
one macroscopic configuration to another. A comparison with other multigrid
sampling schemes (see e.g. \cite{brandt}) reveals that we have derived a reasonably standard scheme, with however
a particularly effective way to store conditional expectations. For details see \cite{chorin9}.

An alternative method for obtaining the expansion coefficients for the renormalized Hamiltonians was proposed in \cite{stinis2}. The method is based on the maximization of the likelihood of the renormalized density. The maximization of the likelihood leads to a moment-matching problem. The moments in this case are the expectation values of the "elementary Hamiltonians" (see above) with respect to the renormalized density. The solution of the moment matching problem yields the expansion of the renormalized Hamiltonian.

The recognition of the links of probability with renormalization is largely due to Jona-Lasinio (see e.g. \cite{jl}). 
The connection of renormalization with incomplete similarity is too well known (see \cite{barenblatt, kadanoff, goldenfeld1})
to require further comment here.

\section{An example: The Korteveg-deVries-Burgers equation}\label{kdv}

As an illustration of the ideas in the previous section, consider
the equation
\begin{equation}
u_t+uu_x=\epsilon u_{xx}-\beta u_{xxx},
\end{equation}
with boundary conditions
\begin{equation}
u(-\infty)=u_0,\ \ u(+\infty)=0, \ \ u_{x}(-\infty)=0,
\end{equation}
where the subscripts denote differentiation, $x$ is the spatial variable,
$t$ is time, $\epsilon>0$ is a diffusion coefficient, $\beta>0$ is a dispersion coefficient and $u_0>0$ is a given constant.
The boundary conditions create a traveling wave solution moving to the right
(towards $+\infty$) with velocity
$u_0/2$ which becomes steady in a moving framework as $t\rightarrow\infty$.
In nondimensional form the equation can be written as: 
\begin{equation}
u_t+uu_x=\frac{1}{R} u_{xx}+u_{xxx},
\label{kdvb}
\end{equation}
with $u_x(-\infty)=0$, $u(+\infty)=0$, $u(-\infty)=1$;
$R=\eps\sqrt U/\alpha$ is a ``Reynolds number".
For $R\leq1$ the traveling wave has a monotonic profile,
while for $R>1$ the profile
is oscillatory, with oscillations whose wave length is of order 1 \cite{bona}.
At zero diffusion $(R=\infty)$
the stationary asymptotic wave train extends to infinity 
on the left. For finite $R$ the wave train is damped and the solution 
tends to 1 as $x$ decreases.

The steady wave profile can be found by noting that it satisfies an ordinary
differential equation, whose solution connects a spiral singularity at $x=\infty$
to a saddle point at $x=+\infty$. 
At the steady state we average the solution at each point $x$ over the region
$\left(x-\ell/2, x+\ell /2\right)$ and call the result $\bar u$. 
Now look for 
an effective equation $g(v,v_x,v_{xx},\ldots)=0$
whose solution $v$ approximates $\bar u$; $v$ can be expected
to be smoother than the solution of (\ref{kdvb}) and thus require fewer mesh points
for an accurate numerical solution.

We now make an analogy between the conditional expectations which define the 
renormalized variables in the previous sections 
and an 
averaging in space which defines ``renormalized"
variables for solutions of the KdVB equations that are stationary
in a moving  frame.
Averaging over an increasing length scale corresponds either to more 
renormalization steps or, equivalently, to renormalization with a greater
number of variables grouped together.
We pick a class of equations in which to seek the ``effective" equation,
the one whose solutions best approximate the averages of the true solution in the
mean square sense; the choice of mean-square approximation
in the KdVB case corresponds to the use of $L_2$ norms implied by the use
of conditional
expectations in the previous sections, and the choice of a class of equations in which to
look for the effective equation is analogous to the choice of a basis
for the representation of the Hamiltonian; the calculation of 
the best coefficients in the chosen class of ``effective" equations corresponds to the
evaluation of the coefficients in the series for the renormalized Hamiltonians. 
In the Hamiltonian case we average the right-hand-sides of the equations and
in the analogous KdVB case we attempt to average the solutions;
this must be so because in the KdVB case we do not have theorems which
guarantee that averaging the right-hand-sides produces the correct statistics for the
solutions.

We can look for an effective equation in the class of equations of
the form
\begin{equation}
-cv_x+vv_x=\epsilon_{eff} v_{xx}+v_{xxx}+\beta |v_x|^\alpha v_{xx}+\dots,
\end{equation}
where $\epsilon\geq0,\alpha\geq0, \beta\geq0$ are constants and $c=1/2$ is the velocity of propagation 
of the steady wave (see also \cite{barenblatt3}).
The problem is to find the value of the parameters  in the effective equation which minimizes
\begin{equation}
I= \int_{-\infty}^{+\infty}|\bar u(x)-v(x)|^2 dx.
\label{min}
\end{equation}
One finds numerically that that the last terms have little effect on the minimum if $I$ when $\ell\ge5$
(in the physics terminology, 
they are ``irrelevant").
The effective equation is thus a Burgers equation 
with a value of the dimensionless diffusion coefficient $\epsilon_{eff}$  different from $1/R$.

The minimization in (\ref{min}) was carried out in \cite{chorin10}, and it showed that the mimimun
was achieved when $\epsilon_{eff}=R^{\nu}\Phi(\ell)$, with the exponent $\nu\sim 0.75$. Note that
when the diffusion coefficient $\epsilon\rightarrow0$,
then $\epsilon_{eff}\rightarrow \infty !$.
This is an incomplete similarity relation, as advertised, relating a ``bare" Reynolds number $R$ to
a ``dressed" Reynolds $\epsilon_{eff}^{-1}$. 
The form of the effective equation could conceivably have been found by averaging the original 
equation, but the relation between the original $\epsilon$ and $\epsilon_{eff}$ requires
some form of renormalization-like reasoning.

\section{The Mori-Zwanzig formalism}\label{mz}
We now return to the problem we started investigating in Section \ref{ave}: How to determine the evolution of
a subset $\hat\varphi$ of components of a vector $\varphi$ described by a nonlinear set of equations
of the form (\ref{eq:system}). This is a nonlinear closure problem of a type much studied in
physics, and a variety of formalisms is available for the job. We choose the Mori-Zwanzig formalism of
irreversible statistical mechanics \cite{fick,grabert,mori,zwanzig,zwanzig2}, because it homes in on the basic difficulty, which is the
description of the memory in the system; the relation of this formalism to other nonlinear formalisms
is described in \cite{CHK04}. That a reduced description of a nonlinear system involves a memory
should be intuitively obvious: suppose you have $n>3$ billiard balls moving about on top of a table
and are trying to describe the motion of just three; the second ball may strike the seventh ball
at a time $t_1$ and the seventh ball may then strike the third ball at a later time. 
The third ball then ``remembers" the state of the system at time $t_1$, and if this memory is
not encoded in the explicit knowledge of where the seventh ball is at all times, then it has to be encoded in some
other way.  We are no longer assuming that the system is Hamiltonian nor that we know an invariant
density.

It is much easier to work with linear equations, and we start by finding a linear equation
equivalent to (not approximating!) the system (\ref{eq:system}). 
Introduce the linear Liouville operator 
$L= \sum_{i=1}^n R_i(x)
\frac{\partial}{\partial x_i}$, and the Liouville equation:
\begin{eqnarray}
\frac{\partial}{\partial t}u(x,t)& = &Lu(x,t) \nonumber\\
u(x,0)& = &g(x),
\label{Liouville}
\end{eqnarray}
with initial data $g(x)$. This is the partial differential
equation for which (\ref{eq:system}) is the set of characteristic equations. One can
verify that the solution of the Liouville equation is $u(x,t)=g(\varphi(x,t))$ (see e.g \cite{CHK}).  In
particular, if $g(x)=x_i$, the solution is $u(x,t)=\varphi_i(x,t),$ the
i-th component of the solution of (\ref{eq:system}). 
This linear partial differential equation is thus equivalent to 
the nonlinear system (\ref{eq:system}). The linearity of equation (\ref{Liouville}) greatly facilitates
the analysis.

Introduce the semigroup notation $u(x,t)=(e^{t L}g)(x)=g(\varphi(x,t))$,
where $e^{tL}$ is the evolution operator associated with the operator $L$;
therefore $e^{tL}g(x)=g(e^{tL}x)$, and
one can also verify that
$e^{tL}L=Le^{tL}$ (this can be seen to be a change of variables formula).  Equation
(\ref{Liouville}) becomes
\[
\frac{\partial}{\partial t}e^{tL}g = L e^{tL} g = e^{tL} Lg.
\]
We suppose that as before we are given 
the initial values of the
$m$ coordinates $\hatx$, and that the distribution of the remaining $n-m$
coordinates $\tilde{x}$ is the conditional density, $f$
conditioned by $\hatx$, where $f$ is initially given.

We define a projection operator $P$ by $Pg=E[g|\hatx]$. 
The conditioning variables are the initial values of $\hat \varphi$;
in section \ref{ave} the conditioning variables were the values of $\hat\varphi(t)$, which are
unusable here when we do not know the probability density at time $t$. Quantities such
as $P\hat \varphi(t)=E[\hat\varphi(t)|\hatx]$ are by definition the best estimates of
the future values of the variables $\hat\varphi$ given the partial data $\hatx$ and are
often the quantities of greatest interest.

Consider 
a resolved coordinate $\varphi_j(x,t)=e^{tL} x_j$ ($j\le m$), and split its time
derivative, $R_j(\varphi(x,t))=e^{tL} L x_j$ as follows: 
\begin{equation}
\frac {\partial}{\partial t} e^{tL} x_j = e^{tL} L x_j =  e^{tL}\P L x_j + e^{tL} \Q L x_j, 
\label{eq:split}
\end{equation}
where $\Q=I-\P$. Define $ \hat{R}_j(\hatx) = (\P R_j)(\hatx)$; the first
term is $e^{tL}\P L x_j =  \hat{R}(\hat{\varphi}(x,t))$ and is a function of the resolved components only (but it is a function of the whole vector of initial data). 
Note that if $Q$ were zero we would recover something that looks 
like the crude approximation of the previous section; however the conditioning
variables are not the same. We shall see that the term in $Q$ is essential.

We further split the remaining term $e^{tL} \Q L x_j$. This splitting will
bring it into a very useful form: a noise term, and a memory term whose kernel depends
on the correlations of the noise term. The fact that such a splitting is possible
is the essence of ``fluctuation-dissipation" theorems (see e.g \cite{landau}).

Let $w(x,t)=e^{t\Q L}\Q L x_j$, i.e., let 
$w(x,t)$ be a solution of the initial value problem: 
\begin{eqnarray}
\frac{\partial}{\partial t}w(x,t)&=&\Q L w(x,t)\  = \  Lw(x,t)- \P Lw(x,t) \nonumber\\
w(x,0)&=&\Q L x_j.
\label{ortho1}
\end{eqnarray}
If for some function h(x), $Ph=0,$ then $Pe^{t\Q L}h=0$ for all time $t$, i.e., $e^{t\Q
L}$ maps the null space of $\P$ into itself.

The evolution operators $e^{tL}$ and $e^{t\Q L}$ satisfy the Duhamel
relation
\[
e^{t L} = e^{t\Q L} + \int_0^t e^{(t-s) L} \P L e^{s \Q L} \,ds.
\]
Hence,
\begin{equation}
e^{tL} Q L x_j = 
e^{t\Q L} \Q L x_j + \int_0^t e^{(t-s)L} \P L e^{s\Q L} \Q L x_j \,ds.
\label{dyson}
\end{equation}

Collecting terms, we find
\begin{equation}
\frac {\partial}{\partial t} e^{tL} x_j =  e^{tL}\P L x_j + 
\int_0^t e^{(t-s) L} \P L e^{s Q L}Q L x_j \,ds +e^{tQL} Q L x_j
\label{eq:langevin}
\end{equation}

The first term on the right hand side is the
Markovian contribution to $\partial_t \varphi_j(x,t)$---it depends only on
the instantaneous value of the resolved $\hatvp(x,t)$.  The second
term depends on $x$ through the values of $\hatvp(x,s)$ at times $s$
between $0$ and $t$, and embodies a memory---a dependence on the past
values of the resolved variables.  Finally, the third term, which
depends on full knowledge of the initial conditions $x$, lies in the
null space of $\P$ and can be viewed as noise with statistics
determined by the initial conditions.

It is important to see that equation (\ref{eq:langevin}) is an identity. The memory and noise
terms have not been added artificially, their presence is a direct consequence of the original
equations of motion. However tempting it may be to average equations by taking one-time 
averages, the results will in general be wrong; one must add a memory and a noise as well.

If what is desired is $P\hat \varphi(t)$, the conditional expectation of
$\hat \varphi(t)$ given $\hat x$ (the best approximation in the sense of $L_2$ to $\hat\vp$ given the
partial data $\hat x$), then one can  premultiply equation (\ref{eq:langevin}) by P; the noise term
then drops out and we find
\begin{equation}
\frac {\partial}{\partial t}P e^{tL} x_j = P e^{tL}\P L x_j + 
P\int_0^t e^{(t-s) L} \P L e^{s Q L}Q L x_j \,ds 
\label{eq:langevin_pro}
\end{equation}
Even if the system we start with is Hamiltonian, the Langevin
equation (\ref{eq:langevin}) is not;  the memory and the noise allow the system to forget
its initial values and decay to ``thermal equilibrium" as it should (see section \ref{ave}).

We now show that the memory term is a functional of the temporal correlations of the noise. 
To save on writing 
we restrict ourselves to cases where the operator $L$ is skew-symmetric,
i.e, $(Lu,v)=-(u,Lv)$, (remember $(u,v)=E[uv]$). The skew-symmetry holds in particular for 
Hamiltonian systems with canonical data, see \cite{CHK3},\cite{evans}; however, here the the assumption is skew-symmetry 
is only an excuse to reduce the number of symbols, 
not a
return to the Hamiltonian case. Pick an orthonormal basis $\{h_k=h_k(\hat x),k=1,\dots\}$ in 
the range of $P$, which is the space of functions of $\hat x$ 
(for example,
the $h_k$ could be Hermite polynomials in the variables $\hatx$). Any function
$\psi(x,t)$,
can be expanded as  $\psi=\sum_k(\psi(x,t),h_k)h_k(\hatx)$, and in particular,
\begin{equation}
P(LQe^{sQL}QLx_j)=\sum_k(LQe^{sQL}QLx_j,h_k)h_k(\hat x).
\label{expand_fin}
\end{equation}
where a factor $Q$ has been inserted before the exponentials, harmlessly because
the operators that follow it all live in the null space of $P$. 
The memory term now becomes
\begin{eqnarray}
\int_0^te^{(t-s)L}PLe^{sQL}QLx_jds\!\!\!&=\!\!\!&\int_0^t\sum_ke^{(t-s)L}(LQe^{sQL}QLx_j,h_k)h_k(\hat x)ds\nonumber\\
\!\!\!&=&\!\!\!\sum_k\!\!\int_0^t(LQe^{sQL}QLx_j,h_k)h_k(\hat \varphi(t-s))ds;
\label{expand}
\end{eqnarray}
In the last identity we used the fact that the parenthesis is independent of time and therefore
commutes with the time evolution operator $e^{tQL}$, and also the fact that $e^{(t-s)L}h_k(\hatx)=h_k(\hat\varphi(t-s))$ by
definition. 
Now $(LQe^{sQL}QLx_j,h_k(\hatx))=-(e^{sQL}QLx_j,QLh_k(\hatx))$ by the symmetry of $Q$
and the assumed skew-symmetry of $L$; each term on the right hand side of equation 
(\ref{expand}) is the ensemble average of the product of the value of the stochastic process $e^{tQL}QLx_j$ at time $s=t$
with the value of the stochastic process $e^{tQL}QLh_k(\hatx)$ evaluated at time $s=0$, i.e., it
is a temporal correlation. All these stochastic processes are in the range of $Q$ for all $t$,
they are therefore components of the noise.  Remember that by definition $Lx_j=R_j$ (a right-hand side in equations (\ref{eq:system})). $PLx_j$ is then an average of the right-hand side of (\ref{eq:system})
and $QLx_j=R_j-E[R_j|\hat x]$ is the initial fluctuation in that right-hand side.

The first, ``Markovian", term in equations (\ref{eq:langevin}) looks straightforward, but perils lurk there
as well. 
In general $R_j$ in equations (\ref{eq:system}) is nonlinear,
and so is $PLx_j=E[R_j|\hat x]$. $e^{tL}PLx_j$ is a nonlinear
function of the functions $\hat\varphi(t)$ which depends on all the components of $x$, not only on $\hat x$.
Some way of approximating this function must be found. If one looks for conditional expectations, one must
find a way to commute $P$ with a nonlinear function; for a discussion, see \cite{CHK3}. This bullet was dodged in section \ref{ave} when the conditioning variables were chosen to be $\hat\varphi(t)$ which change in time, but it may be hard to dodge here.

  																																	      The task now at hand is to extract something usable from these rather cumbersome formulas. A very detailed presentation of
  																																	      the analysis in this section can be found in \cite{c11}.

  \section{Fluctuation-dissipation theorems}\label{fd}

We have established a relation between kernels in the memory term and the noise (the former is made up of covariances of the latter). This is the mathematical content of what are known as ``fluctuation-dissipation theorems" in physics. However, under some specific restricted circumstances, the relation between noise and memory takes on more intuitively appealing forms, which we now briefly describe.
In physics one often takes a restricted basis in the range of $P$ consisting
of the coordinate functions $x_1,...,x_m$ (the components of $\hat{x}$). The resulting projection
is called there the `` linear projection" as if $P$ as defined above were not linear. 
The use of this projection is appropriate when the amplitude of the functions $\hat\phi(t)$ is small. 
One then has 
$h_k(\hat x)=x_k$ for $k\le m$. 
The correlations in equation (\ref{expand}) are then simply 
the temporal correlations of the noise (not of the full solutions of the system!). This is known as the fluctuation-dissipation theorem of the second kind.

Specialize further to a situation where there is a single resolved variable, say $\phi_1$, so that $m=1$
and $\hat\phi$ has a single component. The Mori-Zwanzig equation becomes:

\begin{equation*}
\frac{\partial}{\partial{t}} e^{tL}x_1=
e^{tL}PLx_1+e^{tQL}QLx_1+
\int_0^t e^{(t-s)L}PLe^{sQL}QLx_1ds,
\end{equation*}
or, 
\begin{multline} 
\label{lmz}
\frac{\partial}{\partial{t}} \phi_1(x,t) =
(Lx_1,x_1)\phi_1(x,t)+e^{tQL}QLx_1 \\\
+\int_0^t(LQe^{sQL}QLx_1,x_1)\phi_1(x,t-s)ds \\\
=(Lx_1,x_1)\phi_1(x,t)+e^{tQL}QLx_1-
\int_0^t (e^{sQL}QLx_1,QLx_1)\phi_1(x,t-s)ds, 
\end{multline}
where we have again inserted a harmless factor $Q$ in front of $e^{QL}$, assumed that
$L$ was skew-symmetric as above, and for the sake of simplicity also assumed $(x_1,x_1)=1$
(if the last statement is not true the formulas can be adjusted appropriately). 
Take the inner  product of equation (\ref{lmz}) with $x_1$, you find: 
\begin{multline}
\label{clmz}
\frac{\partial}{\partial{t}} (\phi_1(x,t),x_1)=(Lx_1,x_1)(\phi_1(x,t),x_1) \\\
+(e^{tQL}QLx_1,x_1)-\int_0^t(e^{sQL}QLx_1,QLx_1)\phi_1(x,t-s)ds \\\
=(Lx_1,x_1)(\phi_1(x,t),x_1)-\int_0^t(e^{sQL}QLx_1,QLx_1)
(\phi_1(x,t-s),x_1)ds,
\end{multline}
because $Pe^{tQL}QLx_1=(e^{tQL}QLx_1,x_1)x_1=0$ 
and hence $(e^{tQL}QLx_1,x_1)=0.$ 
Multiply equation (\ref{clmz}) by $x_1$, and remember that  $P\phi_1(x,t)=(\phi_1(x,t),x_1)x_1.$ You find:
\begin{equation}
\label{plmz}
\frac{\partial}{\partial{t}} P\phi_1(x,t)= (Lx_1,x_1)P\phi_1(x,t)-
\int_0^t (e^{sQL}QLx_1,QLx_1)P\phi_1(x,t-s)ds. 
\end{equation}
You observe that the covariance $(\phi_t(x,t),x_1)$ and the projection of $\phi_1$ on $x_1$
obey the same homogenous linear integral equation. This is the fluctuation-dissipation theorem
of the first kind, which embodies the Onsager principle, according to which spontaneous fluctations
in a system
decay at the same rate as perturbations imposed by external means,  when both are small
(so that the linear projection is adequate).
This reasoning can be extended to cases where there are multiple resolved variables, and this is
usually done with the added simplifying assumption that $(x_i,x_j)=0$ when $i\ne j$. We omit the details.

\section{Very short and very long memory approximations}\label{short}

The approximation we shall examine is some detail is:
\begin{equation}
e^{tQL}\cong e^{tL},
\label{QLeL}
\end{equation}
and we will consider under what conditions  it is reasonable. 
We will find that it is reasonable both when memory is very short and when it is very long. The fact that the same approximation works for two opposite cases is not a paradox. The approximation (\ref{QLeL}) states that the orthogonal dynamics operator is very close to the full dynamics operator. In other words, the orthogonal dynamics, which evolve in a space orthogonal to that of the resolved variables, are insensitive to the coupling between resolved and unresolved variables. This can happen in particular when the orthogonal dynamics are very fast or when the orthogonal dynamics are very slow. The ansatz above should work when there is  an effective decoupling of the equations for the resolved and unresolved variables. This raises the question of what determines the range of the memory. Is it possible to have a reduced model with very short or very long memory, depending on how one coarse-grains  a particular system at hand? In \cite{stinis} evidence was presented that, fo!
 r the Kuramoto-Sivashinsky equation, the range of the memory of a reduced model can vary dramatically, depending on whether all the unstable modes in the system are resolved or not. The construction of a reduced model corresponds to renormalization, and the two extreme cases can be interpreted as two fixed points of a renormalization scheme. In which one a reduced model will end up depends on how one renormalizes. Finally, note that the Duhamel formula can be used for an iterative solution of the orthogonal dynamics equation. The term $e^{tL}$ is the zero-th order term of an iterative solution for $e^{tQL}.$ This construction can be based on the use of Feynman diagrams.

First we examine the case when the memory is short, i.e., when the
various terms in the series (\ref{expand_fin}) vanish for $s$ beyond a small value; see \cite{majda} for 
a different approach to short-memory reduced model construction and \cite{stinis3} for comparison with the present short-memory approximation, as well as \cite{p} and the references therein.

The memory term in the Mori-Zwanzig equations (\ref{eq:langevin}) can be rewritten as 
\begin{equation}
\int_0^t e^{(t-s)L} \P L e^{s\Q L} \Q L x_j \,ds =
\int_0^t e^{(t-s)L} \P L \Q e^{s\Q L} \Q L x_j \,ds,
\end{equation}
where the insertion of the extra $\Q$ is harmless.
Adding and subtracting equal quantities, we find:
\begin{equation}
PLe^{sQL}QLx_j=PLQe^{sL} QLx_j + PLQ (e^{sQL}-e^{sL}) QLx_j;
\end{equation}
a Taylor series yields:
\begin{equation}
e^{sQL}-e^{sL}=I+sQL+\dots-I-sL-\dots=-sPL+O(s^2),
\end{equation}
and therefore, using $QP=0$, we find:
\begin{equation}
\int_0^t e^{(t-s)L} P L e^{sQL} Q L x_j \,ds = 
\int_0^t e^{(t-s)L} P L Q e^{sL} Q L x_j\,ds + O(t^3).
\end{equation}
If $P$ is a finite rank projection then
\begin{equation}
P L e^{sQL} Q L x_j = 
\sum_{k} (Q L e^{sQL} Q  L x_j, h_k) h_k(\hatx).
\end{equation}
where, as before, one can write $(QLe^{sQL}QLx_j, h_k)$ as $-(e^{sQL}QLx_j, QLh_k)$ when $L$ is skew-symmetric. 
If the correlations $(e^{sQL}QLx_j,QLh_k)$  and also the correlations $(e^{sL}QLx_j,QLh_k)$ are significant only
over short times $s$, the approximation (\ref{QLeL}) provides an
acceptable approximation without requiring the solution of the
orthogonal dynamics equation (see \cite{stinis} for an application to the dimensional reduction of 
the Kuramoto-Sivashinsky equation and \cite{barber} for an application to molecular dynamics).

The limiting case of the short-memory approximation is when the correlations are delta functions. There is a large literature on solving 
equations (\ref{eq:langevin}) with the
assumption of delta function memory; usually this is done without explicit mention, as if it
were an obvious property of stochastic systems- an astonishing state of affairs
nearly 40 years after Alder and Wainwright demonstrated the long memory 
in a typical physical system \cite{a1}. All the dynamic (i.e., time-dependent)
renormalization group methods we can find depend on this assumption \cite{hohenberg}, and this remark goes a long way towards
explaining their relative lack of success in applications. We will no longer bother making
detailed comparisons with this dynamic renormalization literature; the point of view here is that reduction on the
basis of equations (\ref{eq:langevin}) is the right kind of renormalization, and anything with added drastic assumptions must be justified by appeal to that right kind.

Nevertheless, there are important circumstances where the very short memory assumption can be justified,
in particular in problems with separation of time scales, where the components of $\tilde\varphi(t)$,
the unresolved variables, vary on much faster scales than the resolved variables (see e.g. \cite{majda},\cite{stinis3}).
One can then set 
\begin{equation}
e^{tQL}QLx_j=A_jw_j'(t),
\label{assume}
\end{equation}
where the prime denotes a derivative, the $w_j(t)$ are independent unit Brownian motions,
and the $A_j$ constants that must be derived from some prior knowledge. 
Assume further that the projection $P$ is well represented by the physicists' ``linear" projection and that the density used to perform the projections is invariant. 
The memory term becomes $-A_j^2\delta(t-s)$, equations (\ref{eq:langevin}) become stochastic ordinary
differential equations of the usual kind. As usual (see e.g. \cite{just}), the
corresponding probability densities can be found via Fokker-Planck formalisms (or Kolmogorov
equations, in mathematicians' language). Everything is easier. There is a big literature on 
these methods which we recoil from surveying.

It is often the case that the quantities of interest are the components of $E[\hat\varphi|\hat x]$, and the corresponding projection $P$ is in general poorly approximated by the ``linear" projection. The formalism above readily extends to more general projections, with more terms in the basis chosen in the range of $P$ (see e.g. \cite{CHK3}), as long as one assumes that the temporal correlations of the new terms are fast decaying functions. Terms that have long correlation
times violate the ansatz (\ref{QLeL}) and can hamper rather than enhance accuracy (see e.g. \cite{stinis}). A way to pick the fast decaying terms in the projection of the memory kernel for problems that exhibit separation of time scales was presented in \cite{stinis3}. We should note here that projections which include higher than linear terms are at the heart of mode-coupling theory (see e.g. \cite{schofield}), which has proved very effective in tackling problems in condensed matter physics.

We examine now the validity of the ansatz $e^{tQL}=e^{tL}$ for cases with slowly decaying memory. Write the memory term in the Mori-Zwanzig equation (\ref{eq:langevin}) as
\begin{align*}
\int_0^t e^{(t-s)L}PLe^{sQL}QLx_jds &=\int_0^t Le^{(t-s)L} 
e^{sQL}QLx_jds \\ 
&-\int_0^t e^{(t-s)L}e^{sQL}QLQLx_jds ,
\end{align*} 
where we have used the commutation of $L$ and $QL$ with $e^{tL}$ and $e^{sQL},$ 
respectively. At this point, make the approximation (\ref{QLeL}), which 
eliminates
the $s$ dependence of both integrands and we have
$$\int_0^t e^{(t-s)L}PLe^{sQL}QLx_jds \cong t e^{tL} PLQLx_j.$$
All that remains of the integration in time is the coefficient $t$. 
One can get rid of the noise term by premultiplying equations (\ref{eq:langevin}) by a projection $\P$, as in equation (\ref{eq:langevin_pro}), and obtain a reduced non-autonomous set of differential equations. This approximation was named the $t$-model in \cite{CHK3} (see \cite{ingerman} for an application to the dimensional reduction of a nonlinear Schr\"odinger equation). Other cases where non-Markovian models can be approximated 
by Markovian equations with time-dependent coefficients can be found in \cite{raz}.

We proceed to examine the order of accuracy of this approximation. We have

\begin{multline*}
\int_0^t e^{(t-s)L}PLe^{sQL}QLx_jds- t e^{tL} PLQLx_j = \\\
\int_0^t [e^{(t-s)L}PLe^{sQL}-e^{tL} PL]QLx_jds.
\end{multline*}
Adding and subtracting equal quantities we find

$$ e^{(t-s)L}PLe^{sQL}=e^{tL}PL+e^{tL}[e^{-sL}PLe^{sQL}-PL],$$
and a Taylor series around $s=0$ gives
\begin{equation}\label{t-mod}
e^{-sL}PLe^{sQL}-PL =(I-sL+\ldots)PL(I+sQL+\ldots)-PL=O(s).
\end{equation}
This implies 
$$\int_0^t e^{(t-s)L}PLe^{sQL}QLx_jds=t e^{tL} PLQLx_j + O(t^2).$$
The $O(t^2)$ error estimate can be put into perspective by examining an alternate derivation of the $t$-model. If we expand the integrand of 
the memory term of the Mori-Zwanzig equation around $s=0$ and retain only
the leading term, we find
\begin{align*}
\int_0^t e^{(t-s)L}PLe^{sQL}QLx_jds &= \int_0^t [e^{tL}PLQLx_j 
+O(s)]ds\\
&=t e^{tL} PLQLx_j +O(t^2).
\end{align*}
If we retain only the leading term, we do not keep any information about
the time evolution of the integrand, which in turn means
no information about the evolution of the resolved component and of the 
coupling to the orthogonal dynamics (through the term 
($(LQe^{sQL}QLx_j,h_k)$). Such a drastic approximation is expected to be appropriate in cases where the memory term integrand is slowly decaying, so that information about its initial value is enough.

As an example, consider again the Hald model whose Hamiltonian is
\begin{equation}
H(\phi) = \frac{1}{2} (\phi_1^2 + \phi_2^2 + \phi_3^2 + \phi_4^2 + \phi_1^2 \phi_3^2).
\end{equation}
The resulting equations of motion are:
\begin{align*}
\frac{d\phi_1}{dt} &= \phi_2 \nonumber \\
\frac{d\phi_2}{dt} &= -\phi_1(1 + \phi_3^2) \nonumber \\
\frac{d\phi_3}{dt} &= \phi_4 \nonumber \\
\frac{d\phi_4}{dt} &= -\phi_3(1 + \phi_1^2).
\end{align*}
Suppose one wants to solve only for $\hat\phi=(\phi_1,\phi_2)$, with initial data 
$\hatx=(x_1,x_2)$. Assume the initial data $x_3,x_4$ are sampled from a canonical
density with temperature $T=1$. A quick calculation yields $E[x_3^2|x_1,x_2]=1/(1+x_1^2)$.
the advance in time described by the multiplication by $e^{tL}$ requires just the
substitution $\hatx\rightarrow\hat\phi$. If one commutes the nonlinear function evaluation and
the conditional averaging, i.e., writes $\P f(\hat\phi)=f(\P\hat\phi)$ ( a ``mean-field
approximation"), and writes furthemore $\Phi(t)=\P\hat\phi=E[\hat\phi|\hatx]$ one finds
$\P e^{tL}PLx_1=\Phi_2,\P e^{tL}PLx_2=-\Phi_1(1+1/(1+\Phi_2^2))$; one can calculate
$\P e^{tL}LQLx_j$ for $j=1,2$ and finally one finds:

\begin{align}
\frac{d}{dt}\Phi_1 &=\Phi_2 \nonumber \\
\frac{d}{dt} \Phi_2 &=-\Phi_1 (1 + \frac{1}{1 + \Phi_1^2}) -
2 t \frac{\Phi_1^2 \Phi_2}{(1 + \Phi_1^2)^2}.
\label{eq:hald_t}
\end{align}

The last term represents the damping due to the loss of predictive power
of partial data; the coefficient of the last term increases in time and one may
worry that this last term eventually overpowers the equations and leads to some
odd behavior. This is not the case. Indeed, one can prove the following. If the system
one starts from, equation (\ref{eq:system}) is Hamiltonian with Hamiltonian $H$, and if the
initial data are sampled from an initial canonical density conditioned by partial data $\hat x$,
and if $\hat H$ is the renormalized Hamiltonian ( in the sense of Section \ref{ave}), then
$(d/dt)\hat H \le0$, showing that the components of $\hat\phi$ decay as they should. 
The proof requires a technical assumption ( that the Hamiltonian $H$ can be written
as the sum of a function of $p$ and a function of $q$, a condition commonly satisfied) and
we omit it (see \cite{CHK3}). The reduced system (\ref{eq:hald_t}) was solved numerically in \cite{CHK3}
with gratifying results.

The $t$-model is the zero-th order term in a Taylor expansion (around $s=0$) of the integrand of the memory term in (\ref{eq:langevin}). However, nothing prevents us from keeping more terms in this expansion. Let $$K(\hat{\varphi}(t-s),s)=e^{(t-s)L}PLe^{sQL}QLx_j$$ and expand $K$ around $s=0$, i.e. $$K(\hat{\varphi}(t-s),s)=K(\hat{\varphi}(t),0)+s\frac{\partial K}{\partial s}|_{s=0}+\frac{1}{2}s^2 \frac{\partial^2 K}{\partial s^2}|_{s=0}+O(s^3).$$ In the case when $P$ is the finite-rank projection and the density used to define the projection is invariant, the derivatives of $K$ at $s=0$ are equal-time (static) correlations. In mode-coupling theory, such expressions are known as sum rules. One can assume a functional form for the memory term integrand around $s=0$, e.g. a Gaussian $a e^{-bs^2},$ and use the derivatives of $K$ at $s=0$ to estimate $a,b$ (see \cite{pomeau} for more on sum rules and mode-coupling theory).

\section{Intermediate-range memory}\label{long}

There are intermediate cases where the memory is sufficiently long-range for the short-memory approximation to break down, yet not so slowly decaying that the $t$-model can give accurate results. At present, it is not known how to deal effectively with such cases. In a series of papers \cite{CHK}-\cite{CHK3} we presented special cases and their solutions. In particular in \cite{CHK3} we presented a detailed analysis of the
Hald system. We showed that the memory decays roughly at the same rate as the solution itself (
this is the general case in the absence of separation of scales). We expanded the various correlation functions at equilibrium (i.e., when there are no resolved variables) in Hermite
polynomials, evaluated the coefficients in the expansions by Monte-Carlo once and for all, and then obtained
a system of integro-differential approximations to equations (\ref{eq:langevin}) which we then solved
in various cases. This is a legitimate procedure which may be useful when the same system of equations has to be
solved repeatedly. 
These calculations do exhibit a salient feature of model reduction in time-dependent problems, which is that its set-up costs are often very high. 
The future remedy, if there is one, will surely lie in a deeper understanding of dynamical renormalization and in particular of the
way memory depends on scale.

\section{Acknowledgements} We would like to thank Prof. G.I. Barenblatt, Prof. O. Hald and Prof. R. Kupferman for many helpful 
discussions and comments. 
This work was supported in part by
the National Science Foundation under Grant DMS 04-32710, and by the Director,
Office of Science, Computational and Technology Research,
U.S.\ Department of Energy under Contract No.\ DE-AC03-76SF000098.

\end{document}